\documentclass[11pt,reqno]{amsart}
\usepackage{amsmath}
\usepackage{amsxtra}
\usepackage{amscd}
\usepackage{amsthm}
\usepackage{amsfonts}
\usepackage{amssymb}
\usepackage{eucal}
\let\goth\mathfrak

\newcommand{\nn}{\nonumber}
\newcommand{\bea}{\begin{eqnarray}}
\newcommand{\ena}{\end{eqnarray}}
\newcommand{\be}{\begin{eqnarray*}}
\newcommand{\en}{\end{eqnarray*}}

\renewcommand{\span}{{\mathop{\rm span}}}

\newcommand{\ch}{\mathop{{\rm ch}}}

\newcommand{\bra}[1]{\langle #1 |}        
\newcommand{\ket}[1]{{| #1 \rangle}}      
\newcommand{\slth}{\widehat{\goth{sl}}_2}
\newcommand{\Z}{{\mathbb Z}}  
\newcommand{\C}{{\mathbb C}}

\newcommand{\La}{\Lambda}

\newtheorem{thm}{Theorem}[section]
\newtheorem{prop}[thm]{Proposition}
\newtheorem{lem}[thm]{Lemma}

\numberwithin{equation}{section}

\begin{document}
\pagestyle{myheadings}
\markboth{Feigin, Jimbo, Miwa and Mukhin}
{Symmetric polynomials}

\title[Jack polynomials]
{A differential ideal of symmetric polynomials
spanned by Jack polynomials at $\beta=-(r-1)/(k+1)$}
\author{B. Feigin, M. Jimbo, T. Miwa and E. Mukhin}
\address{BF: Landau institute for Theoretical Physics, Chernogolovka,
142432, Russia}\email{feigin@feigin.mccme.ru}  
\address{MJ: Graduate School of Mathematical Sciences, University of
Tokyo, Tokyo 153-8914, Japan}\email{jimbomic@ms.u-tokyo.ac.jp}
\address{TM: Division of Mathematics, Graduate School of Science, 
Kyoto University, Kyoto 606-8502 Japan}\email{tetsuji@kusm.kyoto-u.ac.jp}
\address{EM: Department of Mathematics, 
Indiana University-Purdue University-Indianapolis, 
402 N.Blackford St., LD 270, 
Indianapolis, IN 46202}
\email{mukhin@math.iupui.edu}

\date{\today}

\setcounter{footnote}{0}\renewcommand{\thefootnote}{\arabic{footnote}}

\begin{abstract}
For each pair of positive integers $(k,r)$ such that $k+1,r-1$ are coprime, 
we introduce an ideal $I^{(k,r)}_n$ of the ring of 
symmetric polynomials $\C[x_1,\cdots,x_n]^{S_n}$. 
The ideal $I^{(k,r)}_n$ has a basis consisting of 
Jack polynomials with parameter $\beta=-(r-1)/(k+1)$, 
and admits an action of a family of differential operators of Dunkl
type including 
the positive half of the Virasoro algebra. 
The space $I^{(k,2)}_n$ coincides with 
the space of all symmetric polynomials in $n$ variables which 
vanish when $k+1$ variables are set equal. The space $I_n^{(2,r)}$
coincides with the space of correlation functions of an abelian
current of a vertex operator algebra related to Virasoro minimal
series $(3,r+2)$.
\end{abstract}
\maketitle



\setcounter{footnote}{0}
\renewcommand{\thefootnote}{\arabic{footnote})}
\renewcommand{\arraystretch}{1.2}

\setcounter{section}{0}
\setcounter{equation}{0}
\section{Introduction}
Let $k$ be a positive integer. 
Consider the space $F^{(k)}_n$ 
of all symmetric polynomials $P(x_1,\cdots,x_n)$ with the following property:
\bea
P(x_1,\cdots,x_n)=0\quad \mbox{if}\quad x_1=\ldots=x_{k+1}.
\label{k+1}
\ena
These polynomials and their analogs 
originate in the work \cite{FS} 
where they were used to study integrable modules over affine Lie algebras. 

To be specific, let $V$ be the level $k$ vacuum module of $\slth$, 
with highest weight vector $\ket{0}$.  
Let $e_n,f_n,h_n$ be the standard generators of 
$\slth$ such that $g_n\ket{0}=0$  ($n\ge 0$, $g=e,f,h$). 
Let further $W=U(\goth{n})\ket{0}$ be the subspace generated from the highest weight vector
by the abelian Lie subalgebra $\goth{n}=\span_{\C}\{e_j\}_{j\in\Z}$. 
Then the dual of the weight $2n$ component $W_n=\{w\in W\mid h_0w=2nw\}$ 
is isomorphic to $F^{(k)}_n$.  
In fact, $F^{(k)}_n$ is nothing but the space of all correlation functions 
of the current operator $e(x)$
\bea
\bra{\psi}e(x_1)\cdots e(x_n)\ket{0}, \quad e(x)=\sum_{n\in\Z}e_n x^{-n-1},   \label{correl}
\ena
where $\bra{\psi}$ runs over the dual space $W_n^*$.
The condition \eqref{k+1} is a consequence of the relation $e(x)^{k+1}=0$, 
which holds in any integrable module of level $k$.  

This explicit realization of $W_n^*$ was utilized in \cite{FS}
to recover a semi-infinite monomial basis of $V$ obtained in \cite{P},  
and to give a representation theoretical interpretation of fermionic
character formulas.  
We remark that \eqref{correl} as well as the condition \eqref{k+1}   
also appear in the literature on 
ground state wave functions in the quantum Hall systems \cite{MR, RR}.  

In this paper we point out a connection between the space 
$F^{(k)}_n$ and Jack polynomials 
which are famous orthogonal polynomials, related to the
Calogero-Sutherland $N$-body problem, see \cite{Mac}, \cite{St}.

Let $k,r$ be positive integers such that $r\ge 2$ and $k+1,r-1$ are
coprime. We consider the Jack polynomials 
specialized to the value of the coupling constant $\beta=-(r-1)/(k+1)$. 
(For our convention on Jack polynomials, see Section \ref{jack sub}).  
In general, Jack polynomials can have poles 
at a negative rational value of $\beta$. However we 
observe that if a partition 
$\lambda=(\lambda_1,\cdots,\lambda_n)$ of length at most $n$ 
satisfies the condition 
\bea
\lambda_i-\lambda_{i+k}\ge r\qquad (1\le i\le n-k),  
\label{cond}
\ena
then the corresponding Jack polynomial 
$P_\lambda(x_1,\dots,x_n)$ does not have a pole at $\beta=-(r-1)/(k+1)$.  
Let $I^{(k,r)}_n$ denote the linear span of $P_\lambda$ satisfying \eqref{cond}. 

We show that 
$I^{(k,r)}_n$ is an ideal of the ring of symmetric polynomials 
$\C[x_1,\cdots,x_n]^{S_n}$, 
and that 
it is stable under the action of a family of linear differential operators  
including the positive half of the Virasoro algebra 
(Theorem \ref{thm:3.1}). 
Using these properties, 
we find that $I^{(k,2)}_n$ coincides with $F^{(k)}_n$ discussed above. 
In other words, the space $F^{(k)}_n$ has a basis consisting of Jack
polynomials evaluated at $\beta=-1/(k+1)$.

We also show that the 
case $k=2$ has a similar interpretation via correlation functions. 
Namely, we verify that $I^{(2,r)}_n$ coincides with the space of 
correlation functions of an abelian current of 
a vertex operator algebra studied in \cite{FJM}. 

The text is organized as follows. 
In Section \ref{sec:2}, after fixing notation, 
we introduce admissible partitions and examine the regularity of $P_\lambda$.
In Section \ref{sec:3} we prove Theorem \ref{thm:3.1} mentioned above.
In Section \ref{sec:4} we examine the special cases 
$r=2$ and $k=2$, respectively. 

\setcounter{section}{1}
\setcounter{equation}{0}
\section{Regularity of Jack polynomials with negative rational coupling constant}\label{sec:2}
\subsection{Jack polynomials}\label{jack sub}
In this section we collect the main
facts about the Jack polynomials
we use in the paper and fix our notations.
For more details on Jack polynomials, see \cite{Mac,St}. 

Let ${S_n}$ be the symmetric group on $n$ letters and
$\La_n=\C(\beta)[x_1,\cdots,x_n]^{S_n}$ 
the ring of symmetric polynomials  
over the field of rational functions $\C(\beta)$.
Let $\pi_n$ denote the set of partitions of length at most $n$.     
To a partition
$\lambda=(\lambda_1,\cdots,\lambda_n)\in\pi_n$, where
$\lambda_i\geq\lambda_{i+1}\geq 0$, we relate the Young
diagram which is 
a subset $\{(i,j)\mid 1\le i\le n,\, 1\le j\le \lambda_i\}$ of $\Z^2$.
An element $(i,j)$ of the Young diagrams is called a node. 
We write $|\lambda|=\sum_{i=1}^n\lambda_i$. 
By $\lambda'$ we denote the partition conjugate to $\lambda$. 
The dominance ordering $\mu\le\lambda$ is defined as 
$\mu_1+\cdots+\mu_i\le \lambda_1+\cdots+\lambda_i$ ($1\le i\le n$).  
{}For $\lambda\in\pi_n$, let 
$m_\lambda=\sum_{\alpha\in S_n\lambda}x^\alpha$
be the orbit sum, where we use the multi-index notation 
$x^\alpha=x_1^{\alpha_1}\cdots x_n^{\alpha_n}$ for 
$\alpha=(\alpha_1,\cdots,\alpha_n)$. 

The Dunkl operators are defined by
\bea
\nabla_i=\partial_i+\beta\sum_{j(\neq i)}\frac{1}{x_i-x_j}(1-K_{ij})
\qquad ( 1\le i\le n).
\label{Dunkl}
\ena
Here 
\bea
(K_{ij}f)(\cdots,x_i,\cdots,x_j,\cdots)=f(\cdots,x_j,\cdots,x_i,\cdots),
\qquad \partial_i=\frac{\partial}{\partial x_i}\notag
\ena 
are the operator which exchanges the $i$-th and the $j$-th variables
and the operator of differentiation with respect to the $i$-th
variable respectively.

We have the commutation relations:
\bea
&&K_{ij}\nabla_m=\nabla_{s_{ij}(m)} K_{ij},\quad 
K_{ij}x_m=x_{s_{ij}(m)} K_{ij},
\nn\\
&&[\nabla_i,\nabla_j]=0, \quad [x_i,x_j]=0, 
\nn\\
&&[\nabla_i,x_j]=\delta_{ij}(1+\beta\sum_{t=1}^nK_{it})-\beta K_{ij}
\label{nabla}
\ena
where $s_{ij}$ signifies the transposition $(i,j)\in S_n$.
Let $D_i=x_i\nabla_i$. The Cherednik operators $\hat{D}_i$ $(1\leq
i\leq n)$ are defined by
\be
\hat{D}_i=D_i+\beta\sum_{j=i+1}^nK_{ij}=x_i\partial_i+\beta\sum_{j(\neq
  i)}\frac{x_{\max \{i,j\}}}{x_i-x_j}(1-K_{ij})+\beta(n-i).
\en
The Cherednik operators also commute with each other 
\be
[\hat{D}_i,\hat{D}_j]=0. 
\en

Define the Sekiguchi operator by 
\bea
S(u,\beta)=\prod_{i=1}^n (u+\hat{D}_i).
\label{Seki}
\ena
In other words, the Sekiguchi operator is the generating function of
the elementary symmetric polynomials in Cherednik operators.

For operators $A$ and $B$, we write $A\sim B$ if $AP=BP$ for any
symmetric function $P$.  For example $\nabla_i\sim \partial_i$,
$D_i\sim\hat{D}_i\sim x_i\partial_i$.

We have 
\be
S(u,\beta)\sim \prod_{i<j}(x_i-x_j)^{-1}\;\det\left[x_i^{n-j}\left(x_i\frac{\partial}{\partial
      x_i}+(n-j)\beta+u\right)\right]_{1\leq i,j \leq n},
\en
see formula (2.21) in \cite{KN} and VI, $\S$ 3, Example 3c in \cite{Mac}.

The action of \eqref{Seki} is triangular on 
$\{m_\lambda\}_{\lambda\in\pi_n}$, 
\bea
S(u,\beta)m_\lambda=
\sum_{\mu\le\lambda}c_{\lambda\mu}(u,\beta)m_\mu, 
\label{triang1}
\ena
where $c_{\lambda\mu}(u,\beta)$ is a polynomial in $u,\beta$. 
In particular, $c_{\lambda\lambda}(u,\beta)$ is given by 
\bea
c_{\lambda\lambda}(u,\beta)=\prod_{i=1}^n(u+\lambda_i+(n-i)\beta).
\label{eigv1}
\ena
The Jack polynomials $\{P_\lambda\}_{\lambda\in\pi_n}$ are the unique 
$\C(\beta)$-basis of $\La_n$ with the following properties:
\bea
&&S(u,\beta)P_\lambda=c_{\lambda\lambda}(u,\beta)P_\lambda, 
\label{eigen}
\\
&&
P_\lambda=
m_\lambda+\sum_{\mu<\lambda}u_{\lambda\mu}(\beta)m_\mu
\qquad (u_{\lambda\mu}(\beta)\in\C(\beta)),
\label{triang2}
\ena
see VI, $\S$ 4, Example 2 of \cite{Mac}.
When necessary we write $P_\lambda(x;\beta)$ for $P_\lambda$, 
exhibiting the $\beta$-dependence explicitly. 

Let 
\bea\label{renorm}
c_\lambda(\beta)=\prod_{(i,j)\in\lambda}((\lambda_j'-i+1)\beta+\lambda_i-j), 
\ena
where the product is over all nodes of partition $\lambda$.
Then the coefficients of the polynomial $J_\lambda=c_\lambda P_\lambda$
are polynomials in $\beta$, see Theorem 3.2 in \cite{KN}. In particular,
the coefficients $u_{\lambda\mu}(\beta)$ are free from poles, except possibly 
at non-positive rational values of $\beta$. 
In Section \ref{reg sec} we describe the properties of $\lambda$ which
are sufficient for regularity of $P_\lambda(x;\beta)$ at a negative
rational value of $\beta$.

The Calogero-Sutherland Hamiltonian is given by
\bea
H=\sum_{i=1}^n(x_i\partial_i)^2+\beta\sum_{i<j}\frac{x_i+x_j}{x_i-x_j}
(x_i\partial_i-x_j\partial_j).
\label{CS}
\ena

We have
\be
H\sim\sum_{i=1}^nD_i^2\sim \sum_{i=1}^n(\hat{D}_i^2-(n-1)\beta\hat{D}_i)+\frac16n(n-1)(n-2)\beta^2.
\en
Therefore we have
\bea\label{CSJ}
HP_\lambda=\varepsilon_\lambda P_\lambda, \qquad
\varepsilon_\lambda=\sum_{i=1}^n(\lambda_i+\beta(n+1-2i))\lambda_i.
\ena

\subsection{Admissible partitions}
In this section we introduce the main combinatorial object of the
paper, the set of admissible partitions.

We call a partition $\lambda\in\pi_n$ {\it $(k,r,n)$-admissible} if 
\bea
\lambda_i-\lambda_{i+k}\ge r
\qquad (1\le i\le n-k).
\label{adm}
\ena
In particular, \eqref{adm} implies 
\be
\lambda_{i}-\lambda_{j}
\ge 
\left[\frac{j-i}{k}\right]r
\en
for all $i<j$, where $[x]$ denotes the largest integer not exceeding $x$. 
 
We will need the following combinatorial lemmas about admissible partitions.
\begin{lem}\label{aux}
Let $\lambda$ be any admissible partition. Then
\bea
(j-i)\beta(k,r)+\lambda_i-\lambda_j & \neq & 0 \qquad (1\leq i<j\le n),\label{den1} \\
(\lambda_j'-i+1)\beta(k,r)+\lambda_i-j & \neq & 0 \qquad (1\leq i\le
n,\quad 1\leq j\leq \lambda_i).\label{den2}
\ena
\end{lem}
\begin{proof}
To get \eqref{den1}, suppose $(j-i)\beta(k,r)+\lambda_i-\lambda_j=0$. 
Then $j-i=(k+1)s$ and $\lambda_i-\lambda_j=(r-1)s$ for some integer $s$. 
Since $j>i$, we get $s>0$.
We have a contradiction:
\be
(r-1)s=\lambda_i-\lambda_j\ge\left[\frac{j-i}{k}\right]r=\left[\frac{(k+1)s}{k}\right]r\ge sr.
\en

To get \eqref{den2}, we suppose $\lambda_j'-i+1=(k+1)s$ and $\lambda_i-j=(r-1)s$ for
some integer $s$. Since $\lambda_i\geq j$, we have $s\geq 0$. The case
$s=0$ means $j=\lambda_i$, $\lambda_j'\geq i$ and therefore is impossible.
For $s>0$, we have a contradiction:
\be
(r-1)s=\lambda_i-j\geq \lambda_i - \lambda_{\lambda_j'} \ge
\left[\frac{\lambda_j'-i}{k}\right]r=\left[\frac{(k+1)s-1}{k}\right]r\ge
sr.
\en
\end{proof}

\begin{lem}\label{aux2}
Let $\lambda$ be $(k,r,n)$-admissible. If $\lambda_j<\lambda_{j-1}$ then
\be
(j-i)\beta(k,r)+\lambda_i-\lambda_j\neq 1 \qquad (i<j). 
\en
\end{lem}
\begin{proof}
Suppose $(j-i)\beta(k,r)+\lambda_i-\lambda_j=1$. 
Then $j-i=(k+1)s$ and $\lambda_i-\lambda_j-1=(r-1)s$ for some integer $s$. 
Since $j>i$, we get $s>0$.
We have a contradiction:
\be
(r-1)s=\lambda_i-\lambda_j-1\ge\lambda_i-\lambda_{j-1}\geq 
\left[\frac{j-i-1}{k}\right]r=\left[\frac{(k+1)s-1}{k}\right]r\ge sr.
\en
\end{proof}

\subsection{Specialization of $P_\lambda$}\label{reg sec}
In this section we examine the regularity of 
Jack polynomials when $\beta$ is negative rational. 

{}Fix a negative rational noninteger number and write it in the form
\bea
-\frac{r-1}{k+1}=:\beta(k,r),
\label{beta0}
\ena
where  $k,r$ are positive integers such that  
$k+1$ and $r-1$ are coprime and $r\ge 2$. 

\begin{lem}\label{adm no pole}
If $\lambda$ is $(k,r,n)$-admissible then $P_\lambda$ has no pole at
$\beta=\beta(k,r)$.
\end{lem}
\begin{proof}
By Lemma \ref{aux}, we have $c_\lambda(\beta(k,r))\neq 0$, where
$c_\lambda$ is given by \eqref{renorm}. The lemma follows. 
\end{proof}

The condition $c_\lambda(\beta(k,r))\neq 0$ is sufficient for
$P_\lambda$ being regular at $\beta=\beta(k,r)$, but not
necessary. The point is the following. Presumably, if the number of
variables in $P_\lambda$ is sufficiently large, then the order of the
pole of $P_\lambda$ in $\beta$ is given exactly by the order of zero
of $c_\lambda$.
However, we have
\be
P_\lambda(x_1,\dots,x_n,0)=P_\lambda(x_1,\dots,x_n),
\en
and in some cases the order of the pole of
$P_\lambda(x_1,\dots,x_n)$ is smaller than that of 
$P_\lambda(x_1,\dots,x_n,x_{n+1})$ and therefore smaller than
the order of zero of $c_\lambda$. This is the case we deal with 
in Proposition \ref{prop:2.1} below. To
prove the regularity of $P_\lambda$ in such a situation we use a
different method.

Before proving the main result of this section, Proposition
\ref{prop:2.1}, we establish a couple of technical lemmas.

\begin{lem}\label{lem:2.1}
Suppose $P_\lambda$ has a pole at $\beta=\beta_0$. 
Then there exists a partition $\nu<\lambda$ such that
\bea
c_{\lambda\lambda}(u,\beta_0)=c_{\nu\nu}(u,\beta_0). 
\label{reson}
\ena
\end{lem}
\begin{proof}
Substituting \eqref{triang1}, \eqref{triang2} into  
\eqref{eigen} and equating coefficients of $m_\nu$, we obtain 
\be
(c_{\nu\nu}(u,\beta)-c_{\lambda\lambda}(u,\beta))u_{\lambda\nu}(\beta)
+\sum_{\nu<\mu<\lambda}u_{\lambda\mu}(\beta)c_{\mu\nu}(u,\beta)
+c_{\lambda\nu}(u,\beta)=0
\en
for all $\nu<\lambda$. 
{}From the assumption, the set of $\mu$ for which 
$u_{\lambda\mu}(\beta)$ has a pole is non-empty.  
A maximal element $\nu$ in this set has the required property \eqref{reson}. 
\end{proof}

\begin{lem}\label{lem:2.2} 
If $P_\lambda$ has a pole at $\beta=\beta(k,r)$, then  
there exists a partition $\nu<\lambda$ and a permutation $w\in S_n$, $w\neq 1$, 
with the properties
\bea
&&
\nu_i=\lambda_{w(i)}+(w(i)-i)\frac{r-1}{k+1},
\qquad (1\le i\le n), 
\label{nu1}
\\
&&w(i)\equiv i~~\bmod~~k+1\qquad (1\le i\le n).
\label{nu2}
\ena
\end{lem}
\begin{proof}
In view of the formula \eqref{eigv1} for $c_{\lambda\lambda}(u,\beta)$, 
the first assertion is an immediate consequence of the previous lemma. 
Since $\nu_i-\lambda_{w(i)}\in\Z$, the second assertion follows.  
\end{proof}

\begin{prop}\label{prop:2.1} Let partition 
$\lambda$ be obtained from a $(k,r,n)$-admissible partition $\mu$   
either by adding or by removing one node. Then
$P_\lambda$ has no pole at $\beta=\beta(k,r)$. 
\end{prop}
\begin{proof}
We have
\bea
\lambda_j-\lambda_{j'}\ge \mu_j-\mu_{j'}-1
\label{ineq1}
\ena
for all $j<j'$. 
Suppose $P_\lambda$ has a pole, 
and take $\nu$ and $w\neq 1$ as in Lemma \ref{lem:2.2}. 
We claim that if $i$ satisfies $w(i)>w(i+1)$, then  
\bea
w(i)=w(i+1)+k, \quad 
\lambda_{w(i+1)}-\lambda_{w(i)}=\mu_{w(i+1)}-\mu_{w(i)}-1=r-1. 
\label{concl}
\ena
  
Indeed, setting $m=w(i)-w(i+1)>0$ we have $m\equiv k\bmod k+1$. 
{}From \eqref{nu1} we obtain 
\bea
\frac{m+1}{k+1}(r-1)\ge
\lambda_{w(i+1)}-\lambda_{w(i)}.
\label{ineq2}
\ena
Since $\mu$ is admissible, we also have 
\bea
\mu_{w(i+1)}-\mu_{w(i)}\ge 
\left[\frac{m}{k}\right]r.
\label{ineq3}
\ena
It follows from \eqref{ineq1}, \eqref{ineq2} and \eqref{ineq3} that 
\be
\frac{m+1}{k+1}(r-1)\ge \left[\frac{m}{k}\right]r-1, 
\en
which is possible only when $m=k$ and \eqref{concl} holds.  

From the assumption,
there is one and only one $i$ 
which violates the condition $w(i)<w(i+1)$.  
We have then $w(j)\ge j$ for $1\le j\le i$.  
Moreover, 
since $w(i)\equiv i\bmod k+1$ and $w(i+1)=w(i)-k$
we have 
$w(i)\ge i+k+1$, $w(i+1)\ge i+1$. 
Since $w(j)<w(j+1)$ holds for $j\ge i+1$ we have also $w(j)\ge j$ for $j\ge i+1$. 
Therefore $w(j)\ge j$ for all $j$. This is a contradiction. 
\end{proof}

Note that if $\lambda$ is as in Proposition \ref{prop:2.1} and if
$\lambda$ is not $(k,r,n)$-admissible then  $c_\lambda$ has a zero of
order one at $\beta=\beta(k,r)$. 

Also note that the proof of Proposition \ref{prop:2.1} with $\mu=\lambda$
gives an alternative proof of Lemma \ref{adm no pole}.

\section{Ideal $I^{(k,r)}_n$ and its properties}\label{sec:3}
\subsection{The ideal $I^{(k,r)}_n$}
In this section we introduce our main object of the study, the space 
$I^{(k,r)}_n$ and describe its properties.

Recall Lemma \ref{adm no pole} which states that,  
when $\lambda$ is a $(k,r,n)$-admissible partition,  
the specialization $P_\lambda$ to $\beta=\beta(k,r)$ given in 
\eqref{beta0} is well-defined as an element of $\C[x_1,\cdots,x_n]^{S_n}$. 
Clearly these polynomials are linearly independent. 
Let $I^{(k,r)}_n$ be their $\C$-linear span, 
\bea
I^{(k,r)}_n=\span_{\C}\{P_\lambda(x_1,\dots,x_n;\beta(k,r))
\mid 
\mbox{$\lambda$ is $(k,r,n)$-admissible}\}.   
\label{I}
\ena

Set 
\be
&&p_m=\sum_{j=1}^n x_j^m \quad (m\ge 1),
\\
&&l_m=\sum_{j=1}^nx_j^{m+1} \partial_j
\qquad (m\ge -1),
\\
&&w^{(t)}_m=\sum_{j=1}^n x_j^{m+t-1} \nabla^{t-1}_j
\qquad (t\ge 2,~~m\ge -t+1),
\en
where $\nabla_j$ are the Dunkl operators \eqref{Dunkl}.

The operators $\{l_m\}_{m\ge -1}$ constitute the positive half of the Virasoro algebra 
\bea
[l_m,l_n]=(n-m)l_{m+n},
\label{vir}
\ena 
and we have 
\be
l_m\sim w^{(2)}_m. 
\en

The operator $w^{(3)}_0$ is related to the Calogero-Sutherland Hamiltonian 
\eqref{CS}
as 
\bea
w^{(3)}_0\sim H+(\beta-1) l_0-\beta p_1 l_{-1}. 
\label{w30}
\ena

\begin{thm}\label{thm:3.1}
\begin{enumerate}
\item $I^{(k,r)}_n$ is an ideal of $\C[x_1,\cdots,x_n]^{S_n}$,  
$p_mI^{(k,r)}_n\subset I^{(k,r)}_n$ ($m\geq 1$).
\item $w^{(t)}_mI^{(k,r)}_n\subset I^{(k,r)}_n$ ($t\ge 2$, $m\ge -t+1$). 
\end{enumerate}
\end{thm}
We defer the proof of Theorem \ref{thm:3.1} 
to Section \ref{proof} and   
mention here an immediate consequence.
\begin{prop}\label{cor:3.1}
Let $P\in I^{(k,r)}_n$. 
Then $(\partial_n^jP)(x_1,\cdots,x_{n-1},0)\in I^{(k,r)}_{n-1}$ 
for all $j\ge 0$.  
\end{prop}
\begin{proof}
Set $(\rho P)(x_1,\cdots,x_{n-1})=P(x_1,\cdots,x_{n-1},0)$. 
Since $\rho(P_\lambda)$ is a Jack polynomial for the same partition $\lambda$ 
in $(n-1)$ variables, 
$\rho$ gives rise to a map 
$I^{(k,r)}_n\rightarrow I^{(k,r)}_{n-1}$. 
The assertion follows from Theorem \ref{thm:3.1} and the relation
\be
\rho\circ \partial^j_n=
\rho\circ\partial_n^{j-1}\circ l_{-1}^{(n)}-l_{-1}^{(n-1)}\circ\rho\circ\partial_n^{j-1},
\en
where we set $l_{-1}^{(n)}=\sum_{j=1}^n\partial_j$.
\end{proof}

\subsection{Proof of Theorem \ref{thm:3.1}}\label{proof}
We prove Theorem \ref{thm:3.1} in several steps. 
First we use 
the following special case of the Pieri formula (\cite{Mac}, VI,(6.7')): 
\bea
&&p_1 P_\mu=\sum_{\lambda}\psi_{\lambda/\mu}'P_\lambda.  
\label{Pieri}
\ena
The sum ranges over partitions 
$\lambda$ obtained by adding one node to $\mu$. 
If $j$ is such that $\lambda_j=\mu_j+1$ holds, then 
$\psi'_{\lambda/\mu}$ is given by 
\bea\label{psi prime}
\psi'_{\lambda/\mu}=\prod_{i=1}^{j-1}
\frac{(j-i-1)\beta+\mu_i-\mu_j}{(j-i)\beta+\mu_i-\mu_j-1}\,
\frac{(j-i+1)\beta+\mu_i-\mu_j-1}
{(j-i)\beta+\mu_i-\mu_j}.
\ena
Equation \eqref{Pieri} is an identity in $\C(\beta)[x_1,\cdots,x_n]^{S_n}$. 

\begin{prop}\label{prop:3.2}
Let $\mu$ be $(k,r,n)$-admissible.   
Then the formula \eqref{Pieri} remains true at $\beta=\beta(k,r)$,  
where we retain only $(k,r,n)$-admissible $\lambda$ in the sum.
In particular 
\be
p_1I^{(k,r)}_n\subset I^{(k,r)}_n.
\en
\end{prop}
\begin{proof}
Let $\lambda$ be a partition appearing in the sum 
\eqref{Pieri}, and let $j$ be such that  
$\lambda_j=\mu_j+1$.  

Note than $P_\lambda$ does not have a
pole at $\beta=\beta(k,r)$ by Proposition \ref{prop:2.1}.

Clearly $\mu_{j-1}>\mu_j$. Then the denominators in \eqref{psi prime}
do not vanish
at $\beta=\beta(k,r)$ by Lemmas \ref{aux} and \ref{aux2}.

Suppose in addition that $\lambda$ is not $(k,r,n)$-admissible.
Then $\mu_{j-k}=\mu_j+r$ and the numerator of the second factor in 
\eqref{psi prime} with $i=j-k$ has a zero at $\beta=\beta(k,r)$.
The proof is over. 
\end{proof}

\begin{prop}\label{prop:3.3} \quad
$
l_{\pm 1}I^{(k,r)}_n\subset I^{(k,r)}_n.
$
\end{prop}
\begin{proof}
We use the following identities in $\C(\beta)[x_1,\dots,x_n]^{S_n}$ 
due to \cite{Las}:
\bea
&&l_1 P_\mu=\sum_{\lambda}\psi''_{\lambda/\mu} P_\lambda,  
\label{Las1}
\\
&&l_{-1}P_\mu=\sum_{\lambda}\widetilde{\psi}'_{\mu/\lambda} P_\lambda.  
\label{Las2}
\ena
In \eqref{Las1} (resp.\eqref{Las2}), 
the sum is taken over $\lambda$ obtained from $\mu$ by adding
(resp. removing) one node. Note that all $P_\lambda$ appearing in 
\eqref{Las1}, \eqref{Las2} have no pole at $\beta=\beta(k,l)$ by
Proposition \ref{prop:2.1}.

In the case \eqref{Las1}, 
\be
\psi_{\lambda/\mu}''=\psi_{\lambda/\mu}'\times 
(-(j-1)\beta+\mu_j)
\en 
where $\lambda_j=\mu_j+1$ and $\psi_{\lambda/\mu}'$ is given by
\eqref{psi prime}. 
Hence the assertion follows from the proof of Proposition \ref{prop:3.2}. 

In the case \eqref{Las2}, 
the formula for $\widetilde{\psi}'_{\mu/\lambda}$ reads 
\bea
\widetilde{\psi}'_{\mu/\lambda}
&=&
\frac1\beta \,((n-i)\beta+\mu_i)
((n-i+1)\beta+\mu_i-1)\times
\nn\\
&&\times 
\prod_{j=i+1}^n
\frac{(j-i-1)\beta+\mu_i-\mu_j}{(j-i)\beta+\mu_i-\mu_j}
\prod_{j=1}^{\mu_i-1}
\frac{(\mu_j'-i+1)\beta+\mu_i-j-1}{(\mu_j'-i+1)\beta+\mu_i-j}, 
\label{Las}
\ena
where $i$ is such that $\lambda_i=\mu_i-1$. 

The denominators in \eqref{Las} do not vanish at $\beta=\beta(k,r)$ by
Lemma \ref{aux}. If in addition $\lambda$ is not $(k,r,n)$-admissible,
then $\mu_{i+k}=\mu_i-r$. 
Since the node $(i,\mu_i)$ is
removable, we have $\mu_i>\mu_{i+1}$. Now the admissibility of $\mu$
forces $\mu_{i+k+1}<\mu_{i+k}$. In particular,
$\mu_{\mu_{i+k}}'=i+k$ and the numerator of the 
factor with $j=\mu_{i+k}$ in the second product of \eqref{Las} vanishes.
\end{proof}

\begin{prop}\label{prop:3.5}\quad
$
l_mI^{(k,r)}_n\subset I^{(k,r)}_n \qquad (m\ge -1).
$
\end{prop}
\begin{proof}
We have $2l_0=[l_{-1},l_1]$. Therefore $l_0I^{(k,r)}_n\subset
I^{(k,r)}_n$ by Proposition \ref{prop:3.3}.

From  \eqref{w30} and \eqref{CSJ} we get $w^{(3)}_0I^{(k,r)}_n\subset
I^{(k,r)}_n$.

Using \eqref{nabla} we find  
\be
&&[w^{(3)}_{0},p_2]\sim 4 l_{2}+2((n-1)\beta+1)p_{2}.
\en
Therefore $l_2I^{(k,r)}_n \subset I^{(k,r)}_n$. 

The proposition then follows from the commutation relations \eqref{vir}. 
\end{proof}

\medskip

\noindent {\it Proof of Theorem \ref{thm:3.1}.}\quad 
Statement (i)
follows from Propositions \ref{prop:3.2}, \ref{prop:3.3},  
and the relation $[l_1,p_m]=mp_{m+1}$ ($m\geq 1$). 

Statement (ii) for $t=2$ follows from Proposition \ref{prop:3.5}.
It remains to show 
\bea
w^{(t)}_mI^{(k,r)}_n\subset I^{(k,r)}_n \quad (t> 2, m\ge -t+1). 
\label{show}
\ena
Using \eqref{nabla} one verifies the following relations. 
\be
&&[l_{-1},w^{(3)}_{0}]=2w^{(3)}_{-1},
\\
&&[w^{(t)}_{-p+1},w^{(3)}_{-1}]=(t-1)w^{(t+1)}_{-t} \qquad (t\geq 2),
\\
&&[w^{(t+1)}_m,p_2]\sim 2tw^{(t)}_{m+2}+t(t-1)(1-\beta)w^{(t-1)}_{m+2}
+2\beta\sum_{i=0}^{t-2}(t-1-i)w^{(i+1)}_{m+t-i}w^{(t-1-i)}_{-t+2+i}.
\en
In the last formula $t\geq 2$ and $m\geq -t$ and we set $w^{(1)}_m=p_m$. 
Since $w^{(2)}_{-1}\sim l_{-1}$, 
the first two formulas imply \eqref{show} for $m+t=1$. 
The general case follows from the third formula by induction on $m+t$. 
\qed

\section{Special cases}\label{sec:4}
In this section we examine the two special cases $r=2$ and $k=2$, 
and identify $I^{(k,r)}_n$ with some spaces of correlation functions. 
For a graded subspace 
$U=\oplus_{d\ge 0}U_d \subset \C[x_1,\cdots,x_n]^{S_n}$, 
the formal character is defined to be
$\ch U=\sum_{d\ge 0}(\dim U_d)q^d$.  
Thus the character of \eqref{I} is 
\be
\ch I^{(k,r)}_n=\sum_{\lambda}q^{|\lambda|}, 
\en
the sum being taken over $(k,r,n)$-admissible partitions $\lambda$.  

\subsection{The case $r=2$}\label{subsec:4.1}
Consider the subspace of symmetric polynomials 
\bea
F^{(k)}_n=\{P\in\C[x_1,\cdots,x_n]^{S_n}
\mid \mbox{$P=0$ if $x_1=\cdots=x_{k+1}$}\}. 
\label{Fkl}
\ena

\begin{prop}\label{prop:3.1} 
$
I^{(k,r)}_n\subset F^{(k)}_n \qquad (n\geq 0).
$
\end{prop}
\begin{proof}
The case $n\le k$ being obvious, we assume that $n\ge k+1$.  

To see the assertion in the case $n=k+1$,  
we use the following specialization formula (see \cite{Mac}, VI,(6.11'))
for Jack polynomials:
\be
P_\lambda(\overbrace{1,\cdots,1}^n;\beta)
=\prod_{(i,j)\in\lambda}
\frac{(n-i+1)\beta+j-1}{(\lambda'_j-i+1)\beta+\lambda_i-j}.
\en
Here the product is taken over all nodes of $\lambda$. 

Let $\lambda$ be $(k,r,n)$-admissible and $n=k+1$.  
Since $\lambda_1-\lambda_{k+1}\ge r$, 
the node $(1,r)$ is contained in $\lambda$.  
This means that the numerator has a zero at $\beta=\beta(k,r)$. 
The denominator does not vanish by Lemma \ref{aux}.
Therefore $P_\lambda=0$ holds 
for $\beta=\beta(k,r)$ and $x_1=\cdots=x_{k+1}$.  

By induction, suppose we have proved the proposition 
for $n-1$, with $n\ge k+2$.
Clearly $P\in F^{(k)}_n$ if and only if 
$(\partial_n^j P)(x_1,\cdots,x_{n-1},0)\in F^{(k)}_{n-1}$ for any $j\ge 0$. 
The assertion then follows from Proposition
\ref{cor:3.1}. 
\end{proof}

\begin{thm}\label{r=2 thm} \quad
$
F^{(k)}_n=I^{(k,2)}_n. 
$
\end{thm}
\begin{proof}
It is known (see \cite{FS}, Proposition 2.6.1' and Theorem 2.7.1) 
that the dimension 
of $(F^{(k)}_n)_d$ is given by 
the number of $(k,2,n)$-admissible partitions such that 
$|\lambda|=d$, and hence $\ch F^{(k)}_n=\ch I^{(k,2)}_n$. 
Now Theorem \ref{r=2 thm} follows from Proposition \ref{prop:3.1}.
\end{proof}

\subsection{The case $k=2$}\label{subsec:4.2} 
Let us consider another special case $k=2$. 
In \cite{FJM}, a vertex operator algebra  
related to the Virasoro minimal series $(3,r+2)$ was studied. 
The main objects in \cite{FJM} are 
an abelian current $a(x)$ which plays a role
analogous to that of $e(x)$ in $\slth$,   
and the `principal subspace' $W$ created 
from the highest weight vector $\ket{0}$ by $a(x)$. 
Denote the space of correlation functions by 
\be
G^{(r)}_n={\rm span}_{\C}
\{\bra{\psi}a(x_1)\cdots a(x_n)\ket{0}
\mid \bra{\psi}\in W^*\}.
\en
The following properties are known:
\begin{enumerate}
\item $\ch G^{(r)}_n=\ch I^{(2,r)}_n$,
\item $G^{(r)}_n$ is generated from the non-zero 
homogeneous component 
of lowest degree by the action of $p_m$ ($m\ge 1$) and $l_m$ ($m\ge -1$), 
\item $P\in G^{(r)}_n$ if and only if 
$(\partial^{\alpha}P)(x_1,x_2,x_3,0,\cdots,0)\in G^{(r)}_3$
holds for all $\alpha=(\alpha_4,\cdots,\alpha_n)\in \Z_{\ge 0}^{n-3}$
where 
$\partial^\alpha=\partial_4^{\alpha_4}\cdots  \partial_n^{\alpha_n}$, 
\end{enumerate}
see Theorem 2.8, Proposition 3.1 and Proposition 4.4 in \cite{FJM}
respectively.

\begin{thm}\label{thm:4.2} \quad
$
G^{(r)}_n=I^{(2,r)}_n.
$
\end{thm}
\begin{proof}
Since both spaces have the same characters, 
it suffices to show the inclusion relation. 

First, consider the case $n=3$.
Let $\varphi_3\in G^{(r)}_3$ be a non-trivial element of the smallest degree $d=r$. 
The corresponding element of $I^{(2,r)}_3$ is 
$\widetilde{\varphi}_3=P_{\lambda}$, with $\lambda=(r,0,0)$. 
For this polynomial we have 
$l_0\widetilde{\varphi}_3=r\widetilde{\varphi}_3$ 
and $H\widetilde{\varphi}_3=\varepsilon_\lambda\widetilde{\varphi}_3$.   
Since removing a node from $\lambda$ leads to non-admissible partitions,  
we have also $l_{-1}\widetilde{\varphi}_3=0$. 
These equations are the same as those known for $\varphi_3$ (see proof
of Proposition 3.2 in \cite{FJM}).  
Since their polynomial solution is unique up to a constant multiple,  
we see that $\varphi_3\in I^{(2,r)}_3$. 
Property (ii) together with Theorem \ref{thm:3.1} 
then imply $G^{(r)}_3\subset I^{(2,r)}_3$, and hence  
$G^{(r)}_3= I^{(2,r)}_3$. 

Now, in the general case $n\ge 4$ the inclusion $I^{(2,r)}_n\subset G^{(r)}_n$
follows from Proposition \ref{cor:3.1} and Property (iii).
This completes the proof. 
\end{proof}

\subsection{Generalizations}
We conclude with two remarks about generalizations of the results of
this paper. 

\medskip

The first remark is about an extension 
to the case of Macdonald polynomials. 

For an indeterminate $s$ we set 
\be
q=s^{k+1}, \quad t=s^{-(r-1)}.
\en
Let 
\be
&&\widetilde{F}^{(k)}_n=
\{P\in\C(s)[x_1,\cdots,x_n]^{S_n}
\mid \mbox{$P=0$ if $x_j=t^{j-1}$ for $j=1,\cdots,k+1$}
\},
\\
&&\widetilde{I}^{(k,r)}_n=\span_{\C(s)}\{P_\lambda(x;q,t)
\mid 
\mbox{$\lambda$ is $(k,r,n)$-admissible}\}.   
\en
Here $P_\lambda(x;q,t)$ denotes the 'monic' 
Macdonald polynomial, see \cite{Mac}, VI,(4.7).

All the working in Subsection \ref{subsec:4.1}
can be extended straightforwardly to get the following theorem.
\begin{thm}\quad
$
\widetilde{F}^{(k)}_n=\widetilde{I}^{(k,2)}_n. 
$
\end{thm}

\medskip

The second remark 
is about the case of general $k,r$. 

It seems natural to anticipate a relation 
similar to Theorem \ref{thm:4.2}. 
Namely we expect that $I^{(k,r)}_n$ coincides with  
the space of correlation functions of 
an abelian current in a vertex operator algebra,  
associated with the minimal series $(k+1,k+r)$ of the $W_k$ algebra. 
We hope to address this subject in the future. 
\bigskip

\noindent
{\it Acknowledgments.}\quad 
J. M. thanks Anatol Kirillov, Masatoshi Noumi, Saburo Kakei and 
Koichi Takemura for helpful information concerning Jack polynomials. 
T. M. thanks Vladimir Fateev and Rinat Kedem for drawing attention to the 
references \cite{MR,RR}. 
This work is partially supported by
the Grant-in-Aid for Scientific Research (B)
no.12440039 and (A1) no.13304010, Japan Society for the Promotion of Science.


\end{document}